\documentclass[11pt]{article}
\usepackage{enumerate}
\usepackage{amssymb,a4wide,latexsym,makeidx,epsfig,fleqn}
\usepackage{amsthm}
\usepackage{amsmath}
\usepackage{enumerate}
\newtheorem{theorem}{Theorem}[section]

\newtheorem{definition}[theorem]{Definition}
\newtheorem{lemma}[theorem]{Lemma}

\begin{document}
\textwidth 150mm \textheight 225mm
\title{The rank of a complex unit gain  graph in terms of the rank of its underlying graph \footnote{This work is supported by the National Natural Science Foundation of China (No. 11171273).}}
\author{{Yong Lu, Ligong Wang\footnote{Corresponding author.} ~and Qiannan Zhou}\\
{\small Department of Applied Mathematics, School of Science, Northwestern
Polytechnical University,}\\ {\small  Xi'an, Shaanxi 710072,
People's Republic
of China.}\\
{\small E-mail: luyong.gougou@163.com, lgwangmath@163.com, qnzhoumath@163.com}}

\date{}
\maketitle
\begin{center}
\begin{minipage}{120mm}
\vskip 0.3cm
\begin{center}
{\small {\bf Abstract}}
\end{center}
{\small Let $\Phi=(G, \varphi)$ be a complex unit gain graph (or $\mathbb{T}$-gain graph)  and $A(\Phi)$ be its adjacency matrix, where $G$ is called the underlying graph of  $\Phi$. The rank of $\Phi$, denoted by $r(\Phi)$, is the rank of $A(\Phi)$. Denote by $\theta(G)=|E(G)|-|V(G)|+\omega(G)$ the dimension of cycle spaces of $G$, where $|E(G)|$, $|V(G)|$ and $\omega(G)$ are the number of edges, the number of  vertices and the number of connected components of $G$, respectively.
In this paper, we investigate bounds for $r(\Phi)$ in terms of $r(G)$, that is, $r(G)-2\theta(G)\leq r(\Phi)\leq r(G)+2\theta(G)$, where $r(G)$ is the rank of $G$. As an application, we also prove that $1-\theta(G)\leq\frac{r(\Phi)}{r(G)}\leq1+\theta(G)$. All corresponding extremal graphs are characterized.

 }
\end{minipage}
\end{center}

\section{Introduction }
In this paper, we only consider simple graphs without multiple edges and loops.
Let $G=(V(G), E(G))$ be a  simple graph with vertex set $V(G)=\{v_{1},v_{2},\ldots,v_{n}\}$ and edge set $E(G)$. The \emph{adjacency matrix} of  $G$  of order $n$ is defined as the $n\times n$ symmetric square matrix $A=A(G)=(a_{ij})_{n\times n}$, where  $a_{ij}=1$  if $v_{i}v_{j}\in E(G)$  and $a_{ij}=0$, otherwise. The \emph{rank} $r(G)$ of $G$ is defined to be the rank of $A(G)$, and the \emph{nullity} $\eta(G)$ of $G$ is defined to be the multiplicity of 0 as an eigenvalue of $A(G)$. Obviously, $r(G)+\eta(G)=n$.
We use Bondy and Murty \cite{BM} for notations  and terminologies not defined here.

Denote by $P_{n},~C_{n}$ a path and a cycle  of order $n$, respectively. A graph is called \emph{empty} if it has some vertices and no edges.
Let $v\in V(G)$, $v$ is called \emph{pendant} vertex if its degree $d_{G}(v)=1$ in $G$, and is called \emph{quasi-pendant} vertex if it is adjacent to a pendant vertex.  An induced subgraph $C_{p}$ of a graph $G$ is called a \emph{pendant cycle} if $C_{p}$ is a cycle and has a unique vertex  of degree 3 in $G$. Denote by $\theta(G)=|E(G)|-|V(G)|+\omega(G)$ the \emph{dimension of cycle spaces} of $G$, where $\omega(G)$ is the  number of connected components of $G$. Obviously, if $\theta(G)=0$, then $G$ is a tree or acyclic. A graph $G$ is called \emph{pairwise vertex-disjoint} if distinct cycles (if any) of $G$ have no common vertices. Let $G$ be a graph with at least one pendant vertex. We call the operation of deleting a pendant vertex and its adjacent vertex from $G$ is a \emph{$\delta$-transformation}. The resultant subgraph $G_{0}$ of $G$ without pendant vertices, obtained from $G$ by applying a series of \emph{$\delta$-transformations}, is called a \emph{crucial subgraph} of $G$.

Let $G$ be a graph with pairwise vertex-disjoint cycles, and let $\mathcal{C}(G)$ denote the set of cycles in $G$. $T_{G}$ is an acyclic graph obtained from $G$ by contracting each cycle $O$ of $G$ into a vertex $t_{O}$. Note that,  $V(T(G))$ is taken to be $U\cup C_{G}$, where $U$ consists of all vertices of $G$ that do not lie on any cycle and $C_{G}=\{t_{O}: O\in \mathcal{C}(G)\}$. Let $u_{1},~u_{2}\in U$, then $u_{1}u_{2}\in E(T_{G})$ if and only if $u_{1}u_{2}\in E(G)$. A vertex $u\in U$ is adjacent to a vertex $t_{O}\in C_{G}$ if and only if $u$ is adjacent (in $G$) to a vertex in the cycle $O$, and two vertices $t_{O_{1}}$, $t_{O_{2}}$ are adjacent in $T_{G}$ if and only if there exists an edge in $G$ joining a vertex of $O_{1}$ to a vertex of $O_{2}$ ($O_{1}, O_{2}\in \mathcal{C}(G)$). Denote by $\Gamma_{G}$ obtained from $T_{G}$ by deleting vertices in $C_{G}$ and all the incident edges.

A complex unit gain graph (or $\mathbb{T}$-\emph{gain graph}) is a graph with the additional structure that each orientation of an edge is given a complex unit, called a \emph{gain}, which is the inverse of the complex unit assigned to the opposite orientation. For a simple graph $G=(V(G),E(G))$ of order $n$, let $\overrightarrow{E}$ be the set of oriented edges, it is obvious that this set contains two copies of each edge with opposite directions. Denote by $e_{ij}$  the oriented edge from $v_{i}$ to $v_{j}$. The \emph{circle group}, which is  denoted by $\mathbb{T}=\{z\in C: |z|=1\}$, is a subgroup of the multiplicative group of all nonzero complex numbers $\mathbb{C}^{\times}$.  A $\mathbb{T}$-gain graph is a triple $\Phi=(G,\mathbb{T},\varphi)$ consisting of a graph $G=(V(G),E(G))$, the circle group $\mathbb{T}$ and a gain function $\varphi:\overrightarrow{E}\rightarrow \mathbb{T}$, where $G$ is the \emph{underlying graph} of $\Phi$  and $\varphi(e_{ij})=\varphi(e_{ji})^{-1}=\overline{\varphi(e_{ji})}$.  We often write $\Phi=(G,\varphi)$ or $G^{\varphi}$ for a $\mathbb{T}$-gain graph. The \emph{adjacency matrix} of the $\mathbb{T}$-gain graph $\Phi$ is the $n\times n$ Hermitian matrix $A(\Phi)=(a_{ij})_{n\times n}$, where $a_{ij}=\varphi(e_{ij})$ if $v_{i}v_{j}\in E(G)$, otherwise $a_{ij}=0$. The \emph{positive inertia index} $i_{+}(\Phi)$, the \emph{negative inertia index} $i_{-}(\Phi)$ and the \emph{nullity} $\eta(\Phi)$ of $\Phi$ are defined to be the number of positive eigenvalues, negative eigenvalues and zero eigenvalues of $A(\Phi)$ including multiplicities, respectively. The \emph{rank} $r(\Phi)$ of a $\mathbb{T}$-gain graph of order $n$ is defined as the rank of $A(\Phi)$. Obviously, $r(\Phi)=n-\eta(\Phi)=i_{+}(\Phi)+i_{-}(\Phi)$. If  every edge of $\Phi$ has gain 1, then $A(\Phi)$ is the same as $A(G)$.

A \emph{subgraph} of $\Phi$ is a subgraph of $G$ in which each edge preserves the original gain in $\Phi$. Let $v\in V(\Phi)$, we write $\Phi-v$ for the induced subgraph obtained from $\Phi$ by deleting the vertex $v$ and all edges incident with $v$. For an induced subgraph $H^{\varphi}$ of $\Phi$, denote by $\Phi-H^{\varphi}$, the subgraph obtained from $\Phi$ by deleting all vertices of $H^{\varphi}$ and all incident edges. For an induced subgraph $F^{\varphi}$ and a vertex $v$ outside $F^{\varphi}$, denote by $F^{\varphi}+v$, the induced subgraph of $\Phi$ with vertex set $V(F)\cup \{v\}$.

Collatz and Sinogowitz \cite{CS} first posed the problem of characterizing nonsingular or singular graph. This problem is of great interest in both chemistry and mathematics. For a bipartite graph $G$ which corresponds to an alternant hydrocarbon in chemistry, if $\eta(G)>0$, it is indicated that the corresponding molecule is unstable. The nullity of a graph is also meaningful in mathematics since it is related to the singularity of adjacency matrix.

For a simple graph, there are  some papers focused on the study of the connections between the nullity (or rank) of graphs $G$ in terms of certain structural parameters, such as matching number, pendant vertices and so on.
Wang and Wong \cite{WW} obtained the  bounds for the nullity of $G$ in terms of the matching number $m(G)$ and $\theta(G)$ of $G$, that is: $|V(G)|-2m(G)-\theta(G)\leq\eta(G)\leq|V(G)|-2m(G)+2\theta(G).$
 Song, Song and Tam \cite{SST} characterized the extreme graphs $G$ that satisfy the upper bound $\eta(G)=|V(G)|-2m(G)+2\theta(G)$ of $G$.  Wang \cite{LWANG} and  Rula et al. \cite{RCZ} independently characterized the lower bound $\eta(G)=|V(G)|-2m(G)-\theta(G)$ of $G$.

In 2016, Ma, Wong and Tian \cite{MWTDAM} have proved the nullity of a graph $G$ in terms of pendent vertices,  that is, $\eta(G)\leq2\theta(G)+p(G),$ where $p(G)$ is the number of pendant vertices of $G$.

For an oriented graph $G^{\sigma}$, there are also some papers studied about the skew-rank of  $G^{\sigma}$.
In 2015, Li and Yu \cite{LXYG} studied the skew-rank of oriented graphs and characterized oriented unicyclic graphs  attaining the minimum value of the skew-rank among oriented unicycle graphs of order $n$ with girth $k$. Qu and Yu \cite{QHYG}, Lu, Wang and  Zhou \cite{LWZ}  characterized the bicyclic oriented graphs with skew-rank 2, 4 and 6, respectively.  Qu, Yu and Feng \cite{QHYGFL} obtained more results about the minimum skew-rank of graphs. They also characterized the unicyclic graphs with skew-rank 4 or 6.

In \cite{CT}, Chen and Tian proved that $sr(G^{\sigma})\geq\sum_{i=1}^{k}q_{i}-2k$ if $G^{\sigma}$ is a connected oriented graph with $k$ pairwise edge-disjoint cycles of size $q_{1},q_{2},\ldots,q_{k}$.  In \cite{MWTLAA}, Ma, Wong and Tian characterized the lower bound and upper bound $2m(G)-2\theta(G)\leq sr(G^{\sigma})\leq2m(G)$ of the skew-rank of an oriented connected graph $G^{\sigma}$ in terms of matching number. The extremal oriented graphs satisfying the lower bound of $sr(G^{\sigma})$ are characterized completely.
In 2016, Wong, Ma and Tian \cite{WMT} have proved that $sr(G^{\sigma})\leq r(G)+2\theta(G)$ for an oriented graph $G^{\sigma}$. They characterized the oriented graphs $G^{\sigma}$ whose skew-rank can attain the upper bound. Lu, Wang and Zhou \cite{LWZ1} characterized the lower bound of the skew-rank of an oriented graphs $G^{\sigma}$, that is $sr(G^{\sigma})\geq r(G)-2\theta(G)$. In 2017, Huang, Li and Wang \cite{HLW} characterized the relation between the skew-rank of an oriented graph and the independence number of its underlying graph.

For a $\mathbb{T}$-gain graph,
Nathan Reff \cite{NR} first defined the adjacency, incidence and Laplacian matrices of a $\mathbb{T}$-gain graph. Yu, Qu and Tu \cite{YQT} characterized some basic properties of positive inertia and negative inertia of a $\mathbb{T}$-gain graph. They also characterized the $\mathbb{T}$-gain unicyclic graphs with small positive or negative index. Lu, Wang and Xiao \cite{LWX} characterized the $\mathbb{T}$-gain bicyclic graphs with rank 2, 3 or 4.

Motivated by the results of relation between the skew-rank of an oriented graph in terms of the rank of its underlying graph, a natural problem is : how about the bounds of the rank of a $\mathbb{T}$-gain graph $\Phi$ in terms of the rank of its underlying graph $G$? In this paper, we will prove that $$r(G)-2\theta(G)\leq r(\Phi)\leq r(G)+2\theta(G).$$  As an application of the lower bound and upper bound of $r(\Phi)$,  we also get that $$1-\theta(G)\leq \frac{r(\Phi)}{r(G)}\leq 1+\theta(G).$$ All corresponding extremal graphs are characterized.

The rest of this paper is organized as follows: In Section 2, we list some known elementary lemmas and results which will be useful in this paper. In Section 3, we characterize the upper bound and lower bound of the rank of a $\mathbb{T}$-gain graph in items of the rank of its underlying graph. In Section 4, we characterize some properties about the upper-optimal and lower-optimal  $\mathbb{T}$-gain graphs, as an application, we also establish sharp upper bound and lower bound of $\frac{r(\Phi)}{r(G)}$.

\section{Preliminaries}

In this section, we introduce some elementary lemmas and known results that will be useful in the future.

 For a $\mathbb{T}$-gain graph $\Phi=(G, \varphi)$, we introduce the following lemmas and results.
\noindent\begin{lemma}\label{le:2.1}\cite{YQT}
\begin{enumerate}[(a)]
  \item Let $\Phi=\Phi_{1}\cup \Phi_{2}\cup \ldots \cup \Phi_{t}$, where $\Phi_{1}, \Phi_{2}, \ldots, \Phi_{t}$ are connected components of a $\mathbb{T}$-gain graph $\Phi$. Then $i_{+}(\Phi)=\sum_{i=1}^{t}i_{+}(\Phi_{i})$.
  \item Let $\Phi$ be a $\mathbb{T}$-gain graph on $n$ vertices. Then $i_{+}(\Phi)=0$ if and only if $\Phi$ is a graph without edges.
  \item Let $H^{\varphi}$ be an induced subgraph of a $\mathbb{T}$-gain graph $G^{\varphi}$. Then $i_{+}(H^{\varphi})\leq i_{+}(G^{\varphi}), i_{-}(H^{\varphi})\leq i_{-}(G^{\varphi})$.
\end{enumerate}
\end{lemma}
\noindent\begin{lemma}\label{le:2.2}\cite{YQT}
Let $\Phi=(G,\varphi)$ be a $\mathbb{T}$-gain graph containing a pendant vertex $v$ with the unique neighbor $u$. Then $i_{+}(\Phi)=i_{+}(\Phi-u-v)+1$, $i_{-}(\Phi)=i_{-}(\Phi-u-v)+1$, $i_{0}(\Phi)=i_{0}(\Phi-u-v)$. Moreover, $r(\Phi)=r(\Phi-u-v)+2$.
\end{lemma}

\noindent\begin{lemma}\label{le:2.3}
If $v$ is a vertex of a $\mathbb{T}$-gain graph $\Phi$, then $r(\Phi)-2\leq r (\Phi-v)\leq r(\Phi)$.
\end{lemma}

\noindent\begin{definition}\label{de:2.4}\cite{LWX}
Let $C_{n}^{\varphi} (n\geq3)$ be a $\mathbb{T}$-gain cycle and $$\varphi(C_{n})=\varphi(v_{1}v_{2}\cdots v_{n}v_{1})=\varphi(v_{1}v_{2})\varphi(v_{2}v_{3})\cdots \varphi(v_{n-1}v_{n})\varphi(v_{n}v_{1}).$$ Then $C_{n}^{\varphi}$ is said to be one of the following five Types:

\begin{displaymath}
\left\{\
        \begin{array}{ll}
          \rm Type~A,&  \emph{if}~\varphi(C_{n})=(-1)^{n/2}~\emph{and}~n~\emph{is~even},\\
          \rm Type~B,& \emph{if}~\varphi(C_{n})\neq(-1)^{n/2}~\emph{and}~n~\emph{is~even},\\
          \rm Type~C,& \emph{if}~Re\left((-1)^{{(n-1)}/{2}}\varphi(C_{n})\right)>0~\emph{and}~n~\emph{is~odd},\\
          \rm Type~D,& \emph{if}~Re\left((-1)^{{(n-1)}/{2}}\varphi(C_{n})\right)<0~\emph{and}~n~\emph{is~odd},\\
          \rm Type~E,& \emph{if}~Re\left((-1)^{{(n-1)}/{2}}\varphi(C_{n})\right)=0~\emph{and}~n~\emph{is~odd},
        \end{array}
      \right.
\end{displaymath}
\end{definition}
where $Re(\cdot)$ is the real part of a complex number.
\noindent\begin{lemma}\label{le:2.5}\cite{YQT}
Let $C_{n}^{\varphi}$ be a $\mathbb{T}$-gain cycle of order $n$. Then

\begin{displaymath}
(i_{+}(C_{n}^{\varphi}),i_{-}(C_{n}^{\varphi}),i_{0}(C_{n}^{\varphi}))=\left\{\
        \begin{array}{ll}
          \left(\displaystyle{\frac{n-2}{2}},\displaystyle{\frac{n-2}{2}},2\right),& \emph{if}~C_{n}^{\varphi}\rm~is~of~Type~A,\\
          \left(\displaystyle{\frac{n}{2}},\displaystyle{\frac{n}{2}},0\right),& \emph{if}~C_{n}^{\varphi}\rm~is~of~Type~B,\\
          \left(\displaystyle{\frac{n+1}{2}},\frac{n-1}{2},0\right),& \emph{if}~C_{n}^{\varphi}\rm~is~of~Type~C,\\
          \left(\displaystyle{\frac{n-1}{2}},\frac{n+1}{2},0\right),& \emph{if}~C_{n}^{\varphi}\rm~is~of~Type~D,\\
          \left(\displaystyle{\frac{n-1}{2}},\frac{n-1}{2},1\right),& \emph{if}~C_{n}^{\varphi}\rm~is~of~Type~E.
        \end{array}
      \right.
\end{displaymath}
\end{lemma}

\noindent\begin{lemma}\label{le:2.6}\cite{YQT}
Let $\Phi=(T, \varphi)$ be a $\mathbb{T}$-gain tree. Then $A(\Phi)$ and $A(T)$ have the same spectrum. Moreover, $r(A(\Phi))=r(A(T))$.
\end{lemma}

 For a simple graph $G$, we introduce the following lemmas and results.

\noindent\begin{lemma}\label{le:2.7}\cite{GX}
Let $G$ be a graph containing a pendant vertex $u$ with the unique neighbor $v$, and $H=G-u-v$ be the induced subgraph of $G$. Then $r(G)=r(H)+2$.
\end{lemma}

\noindent\begin{lemma}\label{le:2.8}\cite{GX}
Let $v$ be a cut-point of a graph $G$ and $G_{1},G_{2},\ldots,G_{t}$ be all  components of $G-v$. If there exists a component, say $G_{1}$, such that $r(G_{1})=r(G_{1}+v)-2$, then $r(G)=r(G-v)+2$. If $r(G_{1})=r(G_{1}+v)$, then $r(G)=r(G_{1})+r(G-G_{1})$.
\end{lemma}


\noindent\begin{lemma}\label{le:2.9}\cite{BBD}
If $v$ is a vertex of a graph $G$, then $r(G)-2\leq r (G-v)\leq r(G)$.
\end{lemma}

\noindent\begin{lemma}\label{le:2.10}\cite{CL}
Let $C_{p}$ be a cycle of order $p$. Then $r(C_{p})=p-2$ if $p\equiv0(\emph{mod}~4)$, and $r(C_{p})=p$ otherwise.
\end{lemma}

Let $T$ be an acyclic graph with at least one edge, we denote by $\widetilde{T}$ the subgraph obtained from $T$ by deleting all pendant vertices of $T$.
\noindent\begin{lemma}\label{le:2.11}\cite{MWTLAA}
 Let $T$ be an acyclic graph with at least one edge. Then
 \begin{enumerate}[(a)]
  \item  $r(\widetilde{T})<r(T)$.
  \item If $r(T-W)=r(T)$ for a subset $W$ of $V(T)$, then there is a pendant vertex $v$ such that $v\notin W$.
\end{enumerate}
\end{lemma}

\section{Relation between the rank of a $\mathbb{T}$-gain graph and the rank of its underlying graph}
In this section, we will give the upper bound and lower bound of $r(\Phi)$ of a $\mathbb{T}$-gain graph $\Phi=(G, \varphi)$  in terms of $r(G)$ and $\theta(G)$. First, we will introduce the following lemma that will be useful for later.
\noindent\begin{lemma}\label{le:3.1}\cite{WMT}
Let $G$ be a  graph with a vertex $v$. Then
\begin{enumerate}[(a)]
  \item  $\theta(G)=\theta(G-v)$ if $v$ lies outside any cycle of $G$.
  \item $\theta(G-v)\leq \theta(G)-1$ if $v$ lies on a cycle of $G$.
  \item $\theta(G-v)\leq \theta(G)-2$ if $v$ is a common vertex of distinct cycles of $G$.
  \item If the cycles of $G$ are pairwise vertex-disjoint, then $\theta(G)$ precisely equals the number of cycles in $G$.
\end{enumerate}
\end{lemma}
From \cite{WMT}, we know that a similar result as Lemma \ref{le:3.1} holds for a $\mathbb{T}$-gain  $\Phi=(G, \varphi)$.

\noindent\begin{theorem}\label{th:3.2}
Let $\Phi=(G, \varphi)$ be a $\mathbb{T}$-gain graph of order $n$. Then $$r(G)-2\theta(G)\leq r(\Phi)\leq r(G)+2\theta(G).$$
\end{theorem}
\noindent\textbf{Proof.}
 We shall apply induction on $\theta(G)$ to prove the lower bound and upper bound of $r(\Phi)$.

\textbf{Case 1.} If $\theta(G)=0$, that is $\Phi$ is a $\mathbb{T}$-gain acyclic, then the result follows from Lemmas \ref{le:2.1} and \ref{le:2.6}.

\textbf{Case 2.} If $\theta(G)\geq1$, i.e., there is at  least one cycle in $\Phi$. Let $v$ be a vertex of a cycle of $\Phi$. By Lemma \ref{le:3.1},
\begin{equation}\label{1}
\theta(G-v)\leq \theta(G)-1.
\end{equation}
The induction hypothesis to $\Phi-v$ means that
\begin{equation}\label{2}
r(G-v)-2\theta(G-v)\leq r(\Phi-v)\leq r(G-v)+2\theta(G-v).
\end{equation}
By Lemmas \ref{le:2.3} and \ref{le:2.9}, we have
$r(\Phi-v)\leq r(\Phi)\leq r(\Phi-v)+2,~r(G-v)\leq r(G)\leq r(G-v)+2.$

So,
\begin{equation}\label{3}
r(\Phi)\geq r(\Phi-v)\geq r(G-v)-2\theta(G-v)\geq r(G)-2-2\theta(G)+2=r(G)-2\theta(G),
\end{equation}
and
\begin{equation}\label{4}
r(\Phi)\leq r(\Phi-v)+2\leq r(G-v)+2\theta(G-v)+2\leq r(G)+2\theta(G)-2+2=r(G)+2\theta(G).
\end{equation}
That is, the lower bound and upper bound of $r(\Phi)$ are both obtained.

This completes the proof. \quad $\square$

For convenience, we call a $\mathbb{T}$-gain graph $\Phi=(G, \varphi)$ \emph{lower-optimal} (\emph{upper-optimal}) if $r(\Phi)$ can obtain the  lower bound (upper bound) in Theorem \ref{th:3.2}.

\section{The necessary and  sufficient conditions of the rank of a $\mathbb{T}$-gain graph $\Phi$ which attains the upper bound and lower bound of $r(\Phi)$}

In this section, we will give some useful lemmas and the necessary  and  sufficient conditions of the rank of a $\mathbb{T}$-gain graph $\Phi=(G, \varphi)$ which attains the upper bound and lower bound of $r(\Phi)$.

\noindent\begin{lemma}\label{le:4.1}
Let $\Phi=(G, \varphi)$ be a $\mathbb{T}$-gain graph and  $v$ be a vertex lying on a cycle of $\Phi$.
\begin{enumerate}[(i)]
  \item  If $\Phi$ is lower-optimal, then $r(\Phi)=r(\Phi-v)$,~$r(G-v)=r(G)-2$,~$\theta(G)=\theta(G-v)+1$.
   \item If $\Phi$ is upper-optimal, then $r(\Phi)=r(\Phi-v)+2$,~$r(G-v)=r(G)$,~$\theta(G)=\theta(G-v)+1$.
  \item If $\Phi$ is lower-optimal  (upper-optimal), then  $\Phi-v$ is lower-optimal (upper-optimal), and $v$ lies on only one cycle of $G$ and $v$ is not a quasi-pendant vertex of $G$.
\end{enumerate}
\end{lemma}

\noindent\textbf{Proof.}
If $\Phi$ is lower-optimal, we have
$r(G)-2\theta(G)\leq r(\Phi)=r(G)-2\theta(G)$, where the inequality follows from Theorem \ref{th:3.2}. So, the inequalities (\ref{1}), (\ref{3}) and the lower bound of inequality (\ref{2}) in the proof of Theorem \ref{th:3.2} all turn into equalities. That is, (i) and the first assertion of (iii) of this lemma are obtained. Furthermore,  combining with $\theta(G-v)=\theta(G)-1$ and (c) of Lemma \ref{le:3.1}, $v$ lies on only one cycle of $G$.

Suppose that $v$ is a quasi-pendant vertex adjacent to a pendant vertex $u$,  we have $r(\Phi)=r(\Phi-v)=r(\Phi-v)+2$, where the first equality follows from (i) of this lemma, and the second equality follows from  Lemma \ref{le:2.2}, a contradiction.

If $\Phi$ is upper-optimal, we have
$r(G)+2\theta(G)=r(\Phi)\leq r(G)+2\theta(G)$, where the inequality follows from Theorem \ref{th:3.2}. So, the inequalities (\ref{1}), (\ref{4}) and the upper bound of inequality (\ref{2})  in the proof of Theorem \ref{th:3.2} all turn into equalities. That is, (ii) of this lemma is obtained and $\Phi-v$ is also upper-optimal. Furthermore, $v$ lies on only one cycle of $G$.

Similar as above, suppose that $v$ is a quasi-pendant vertex of $G$. By Lemma \ref{le:2.7} and (ii) of this lemma, we have $r(G)=r(G-v)=r(G-v)+2$, a contradiction.

This completes the proof. \quad $\square$

\noindent\begin{theorem}\label{th:4.2}
Let $\Phi=(G, \varphi)$ be a $\mathbb{T}$-gain graph with a pendant vertex $u$ adjacent to $v$. Let $F^{\varphi}=\Phi-u-v$. If $\Phi$ is upper-optimal (lower-optimal), then $v$ cannot  lie on any cycle of $G$, and $F^{\varphi}$ is also upper-optimal (lower-optimal).
\end{theorem}
\noindent\textbf{Proof.}
By (iii) of Lemma \ref{le:4.1}, we know that if $v$ is a quasi-pendant vertex of $G$, then $v$ does not lie on any cycle of $G$, i.e., $\theta(G)=\theta(G-v)$.

If $\Phi$ is lower-optimal, by  Lemmas \ref{le:2.2} and \ref{le:2.7}, we have
$$r(\Phi)=r(G)-2\theta(G)=r(G-v)+2-2\theta(G-v)=r(F^{\varphi})+2.$$

That is, $r(F^{\varphi})=r(G-v)-2\theta(G-v)$, i.e., $F^{\varphi}$ is  lower-optimal.

If $\Phi$ is upper-optimal, by  Lemmas \ref{le:2.2} and \ref{le:2.7}, we have
$$r(\Phi)=r(G)+2\theta(G)=r(G-v)+2+2\theta(G-v)=r(F^{\varphi})+2.$$

That is, $r(F^{\varphi})=r(G-v)+2\theta(G-v)$, i.e., $F^{\varphi}$ is  upper-optimal.

This completes the proof. \quad $\square$

In fact, from the process of proof in Theorem \ref{le:2.8} of \cite{LWX}, we can see that when $n$ is odd and $a\neq0$, i.e., $C_{n}^{\varphi}$ is of Type C or D, then $r(\Phi)\geq n-1+r(H)$ (Since $r(C)\geq r(H)$). So, we rewrite  the result of the  this theorem when $C_{n}^{\varphi}$ is of Type C or D.
\noindent\begin{lemma}\label{le:4.3}\cite{LWX}
Let $C_{n}^{\varphi_{0}}$ be a $\mathbb{T}$-gain cycle of order $n~(n\geq3)$ and $H=(G_{1},\varphi_{1})$ be a $\mathbb{T}$-gain graph of order $m~(m\geq1)$. Assume that $\Phi=(G,\varphi)$ is the $\mathbb{T}$-gain graph obtained by identifying a vertex of $C_{n}^{\varphi_{0}}$ with a vertex of $H$ (i.e., $V(C_{n}^{\varphi_{0}})\cap V(H)=v $). Let $F=(G_{2},\varphi_{2})$ be the induced subgraph obtained from $H$ by deleting the vertex $v$ and its incident edges. Then

\begin{displaymath}
\left\{\
        \begin{array}{ll}
          r(\Phi)=n-2+r(H),& \emph{if}~C_{n}^{\varphi}~\emph{is~of~Type~A},\\
          r(\Phi)=n+r(F), & \emph{if}~C_{n}^{\varphi}~\emph{is~of~Type~B},\\
          r(\Phi)=n-1+r(H),& \emph{if}~C_{n}^{\varphi}~\emph{is~of~Type~E},\\
          n-1+r(H)\leq r(\Phi)\leq n+r(H), & \emph{if}~C_{n}^{\varphi}~\emph{is~of~Type~C~or~D}.
        \end{array}
      \right.
\end{displaymath}
\end{lemma}




\noindent\begin{theorem}\label{th:4.4}
Let $\Phi=(G, \varphi)$ be a $\mathbb{T}$-gain graph and $C_{p}^{\varphi}$ be a  $\mathbb{T}$-gain pendant cycle of order $p$, $d(v)=3, v\in V(C_{p})$, and let $F^{\varphi}=(F, \varphi)=\Phi-C_{p}^{\varphi}$, $H^{\varphi}=(H, \varphi)=F^{\varphi}+v$.
If $\Phi$ is lower-optimal, then
\begin{enumerate}[(a)]
  \item   $C_{p}^{\varphi}$ is of Type A with order $p\equiv 2(\emph{mod}~4)$.
  \item  Both $H^{\varphi}$ and $F^{\varphi}$ are lower-optimal.
  \item $r(\Phi)=p-2+r(F^{\varphi})$, $r(F^{\varphi})=r(H^{\varphi})$, $r(G)=p+r(H)$ and $r(F)=r(H)$.
\end{enumerate}
\end{theorem}
\noindent\textbf{Proof.} 
\textbf{Claim 1.} $p$ is even.

Note that $\theta(G)=\theta(H)+1$, suppose that $p$ is odd, then $C_{p}^{\varphi}$ is of Type C, D or E. We have
$r(G)=r(\Phi)+2\theta(G)\geq p-1+r(H^{\varphi})+2\theta(G)\geq p-1+r(H)-2\theta(H)+2\theta(G)=p-1+r(H)+2=p+1+r(H)$, where the first inequality follows from Lemma \ref{le:4.3}, the second inequality is application of the lower bound of Theorem \ref{th:3.2}.

On the other hand, by (i) of Lemma \ref{le:4.1}, we have
\begin{equation}\label{5}
r(G)=r(G-v)+2=p-1+r(F)+2=p+1+r(F).
\end{equation}
So, we have $r(F)\geq r(H)$ (note that $r(G)\geq p+1+r(H)$).

By Lemma \ref{le:2.9}, we know that $r(F)\leq r(H)$. So, $r(F)=r(H).$

Let $A(G)$ be the adjacency matrix of $G$, where
\begin{displaymath}
 A(G)=\left(
  \begin{array}{cccccccccccccc}
            A&     \alpha&   0 \\
            \alpha^{T}&     0&   \beta\\
            0&    \beta^{T}&   B\\

  \end{array}
\right),
\end{displaymath}
where $A$ is the adjacency matrix of $C_{p}-v$, $B$ is adjacency matrix of $F$, $\alpha^{T}$ refers to the transpose of $\alpha$.
From the process of proof in the Lemma 4.4 in \cite{WMT}, we have
 \begin{displaymath}
 r(G)=r\left(
  \begin{array}{cccccccccccccc}
            A&     0&   0 \\
            0&     a&   \beta\\
            0&    \beta^{T}&   B\\

  \end{array}
\right),
\end{displaymath}
where $a=-\alpha^{T}A^{-1}\alpha$. Note that $p$ is odd, and $r(A)=r(C_{p}-v)=r(P_{p-1})=p-1$.
So, \begin{displaymath}
 r(G)=r(A)+r\left(
  \begin{array}{cccccccccccccc}
               a&   \beta\\
               \beta^{T}&   B\\
  \end{array}
\right).
\end{displaymath}
Then we have
 \begin{displaymath}
 r(G)=r(A)+r\left(
  \begin{array}{cccccccccccccc}
               a&   \beta\\
               \beta^{T}&   B\\
  \end{array}
\right)\leq r(A)+r\left(
  \begin{array}{cccccccccccccc}
            a&         0\\
            0&          0 \\
  \end{array}
  \right)+r\left(
  \begin{array}{cccccccccccccc}
             0&   \beta\\
               \beta^{T}&   B\\
  \end{array}
  \right).
\end{displaymath}

That is  $r(G)\leq p-1+1+r(H)=p+r(H)=p+r(F)$, which contradicts to Equation (\ref{5}). So, $p$ is even.

\textbf{Claim 2.} Both $H^{\varphi}$ and $F^{\varphi}$ are lower-optimal.




Let $u$ be a vertex of $C_{p}$ adjacent to $v$. By (i) of Lemma \ref{le:4.1}, Lemmas \ref{le:2.2}, \ref{le:2.7} and the fact that $p$ is even, we have
\begin{equation}\label{6}
 r(\Phi)=r(\Phi-u)=p-2+r(H^{\varphi}).
\end{equation}
\begin{equation}\label{7}
 r(G)=r(G-u)+2=p-2+r(H)+2=p+r(H).
\end{equation}

Since $u$ lies on $C_{p}$, so
\begin{equation}\label{8}
 \theta(G)=\theta(H)+1=\theta(F)+1.
\end{equation}

Combining with Equations (\ref{6})--(\ref{8}), we have
$r(\Phi)=r(G)-2\theta(G)=p+r(H)-2\theta(H)-2=p-2+r(H^{\varphi})$, so
 $r(H^{\varphi})=r(H)-2\theta(H).$
That is, $H^{\varphi}$ is  lower-optimal.

 By (i) of Lemma \ref{le:4.1} and Lemmas \ref{le:2.2} and \ref{le:2.7}, we also have
 \begin{equation}\label{9}
 r(\Phi)=r(\Phi-v)=p-2+r(F^{\varphi}).
\end{equation}
\begin{equation}\label{10}
 r(G)=r(G-v)+2=p-2+r(F)+2=p+r(F).
\end{equation}

Combining with Equations (\ref{7}) and (\ref{10}), we have $r(H)=r(F)$.

Combining with Equations (\ref{8}), (\ref{9}) and (\ref{10}), we have
$r(\Phi)=r(G)-2\theta(G)=p+r(F)-2\theta(F)-2=p-2+r(F^{\varphi})$, so,
 $r(F^{\varphi})=r(F)-2\theta(F)$.
That is, $F^{\varphi}$ is also lower-optimal.


\textbf{Claim 3}. $p\equiv 2(\rm{mod}~4)$.

Suppose to the contrary that $p\equiv 0(\rm{mod}~4)$. Denote by $C_{p}-v=P_{p-1}$, by Lemmas \ref{le:2.7} and \ref{le:2.10}, we have $r(P_{p-1})=r(C_{p})$. By Lemma \ref{le:2.8}, we have
$r(G)=r(P_{p-1})+r(H)=p-2+r(H)$, which contradicts to Equation (\ref{7}).

\textbf{Claim 4}. $C_{p}^{\varphi}$ is of Type A.

Suppose to the contrary that $C_{p}^{\varphi}$ is of Type B, then by Lemma \ref{le:4.3}, we have
 $r(\Phi)=p+r(F^{\varphi})$,
which contradicts to Equation (\ref{9}).

This completes the proof. \quad $\square$


\noindent\begin{theorem}\label{th:4.5}
Let $\Phi=(G, \varphi)$ be a $\mathbb{T}$-gain graph and $C_{p}^{\varphi}$ be a  $\mathbb{T}$-gain pendant cycle of order $p$, $d(v)=3, v\in V(C_{p})$, and let $F^{\varphi}=(F, \varphi)=\Phi-C_{p}^{\varphi}$, $H^{\varphi}=(H, \varphi)=F^{\varphi}+v$.
If $\Phi$ is upper-optimal, then
\begin{enumerate}[(a)]
  \item  $C_{p}^{\varphi}$ is of Type B with order $p\equiv 0(\emph{mod}~4)$.
  \item  Both $H^{\varphi}$ and $F^{\varphi}$ are upper-optimal.
  \item $r(\Phi)=p+r(F^{\varphi})$, $r(F^{\varphi})=r(H^{\varphi})$, $r(G)=p-2+r(H)$ and $r(F)=r(H)$.
\end{enumerate}
\end{theorem}
\noindent\textbf{Proof.} 
\textbf{Claim 1.} $C_{p}^{\varphi}$ can not be of Type C, D or E.

At first, suppose that $C_{p}^{\varphi}$ is of Type C or D, by Lemma \ref{le:4.3} and the application of Theorem \ref{th:3.2} to $H^{\varphi}$, we have

$$r(G)+2\theta(G)=r(\Phi)\leq p+r(H^{\varphi})\leq p+r(H)+2\theta(H).$$

So, we have
$r(G)\leq r(H)+p-2$ (note that $\theta(G)=\theta(H)+1$). Note that $p$ is odd, from (ii) of Lemma \ref{le:4.1}, we have
 \begin{equation}\label{11}
 r(G)=r(G-v)=p-1+r(F).
\end{equation}
 That is, $p-1+r(F)\leq r(H)+p-2$, i.e.,
$r(F)\leq r(H)-1$.

On the other hand, by Lemma \ref{le:2.9}, we have $r(F)\geq r(H)-2$, then $r(H)-2\leq r(F)\leq r(H)-1.$
That is, $r(F)=r(H)-1$ or $r(F)=r(H)-2$.

Let $A(G)$ be the adjacency matrix of $G$ as described in Theorem \ref{th:4.4}.

\textbf{Case 1.} $r(F)=r(H)-1$.

\textbf{Subcase 1.1.} The vector $\beta$ can be linearly expressed by the row vector of $B$, then
\begin{displaymath}
 r(H)=r\left(
  \begin{array}{cccccccccccccc}
             0&   \beta\\
               \beta^{T}&   B\\
  \end{array}
  \right)\leq r\left(
  \begin{array}{cccccccccccccc}
             a&   \beta\\
               \beta^{T}&   B\\
  \end{array}
  \right)\leq r\left(
  \begin{array}{cccccccccccccc}
             a&   0\\
               \beta^{T}&   0\\
  \end{array}
  \right)+ r\left(
  \begin{array}{cccccccccccccc}
             0&   \beta\\
               0&   B\\
  \end{array}
  \right)=1+r(F).
\end{displaymath}

Combining with $r(F)=r(H)-1$, we have
\begin{displaymath}
 r\left(
  \begin{array}{cccccccccccccc}
               a&   \beta\\
               \beta^{T}&   B\\
  \end{array}
\right)=1+r(F).
\end{displaymath}

So, combining with Lemma \ref{le:2.7}, we have $r(G)=r(A)+1+r(F)=p-1+1+r(F)=p+r(F)$,
which contradicts to Equation (\ref{11}).

\textbf{Subcase 1.2.} The vector $\beta$ cannot be linearly expressed by the row vector of $B$,
then
\begin{displaymath}
 r(H)=r\left(
  \begin{array}{cccccccccccccc}
             0&   \beta\\
               \beta^{T}&   B\\
  \end{array}
  \right)\leq r\left(
  \begin{array}{cccccccccccccc}
             a&   \beta\\
               \beta^{T}&   B\\
  \end{array}
  \right)\leq r\left(
  \begin{array}{cccccccccccccc}
             a&   0\\
               0&   0\\
  \end{array}
  \right)+ r\left(
  \begin{array}{cccccccccccccc}
             0&   \beta\\
               \beta^{T}&   B\\
  \end{array}
  \right)=1+r(H).
\end{displaymath}

So, we have
\begin{displaymath}
 r\left(
  \begin{array}{cccccccccccccc}
               a&   \beta\\
               \beta^{T}&   B\\
  \end{array}
\right)=r(H)~\textrm{or}~1+r(H).
\end{displaymath}

Hence, $r(G)=r(A)+r(H)=r(A)+r(F)+1=p+r(F)$, or $r(G)=r(A)+r(H)+1=r(A)+r(F)+2=p+1+r(F)$, both contradict to Equation (\ref{11}).

\textbf{Case 2.} $r(F)=r(H)-2$.

We say the vector $\beta$ cannot be linearly expressed by the row vector of $B$. Otherwise, \begin{displaymath}
 r(F)+2=r(H)=r\left(
  \begin{array}{cccccccccccccc}
             0&   \beta\\
               \beta^{T}&   B\\
  \end{array}
  \right)\leq r\left(
  \begin{array}{cccccccccccccc}
             0&   0\\
               \beta^{T}&   0\\
  \end{array}
  \right)+ r\left(
  \begin{array}{cccccccccccccc}
             0&   \beta\\
               0&   B\\
  \end{array}
  \right)=1+r(F),
\end{displaymath}
a contradiction.
 Then, similar to the proof in Subcase 1.2, we  can also get  a contradiction to Equation (\ref{11}).
So, $C_{p}^{\varphi}$ cannot be  of Type C or D.

If $C_{p}^{\varphi}$ is of Type E, then by Lemma \ref{le:4.3}, we have
$$r(G)+2\theta(G)=r(\Phi)=p-1+r(H^{\varphi})\leq p-1+r(H)+2\theta(H).$$

So, we have $r(G)\leq p-3+r(H)$, combining with Equation (\ref{11}), we have $r(F)+2\leq r(H)$. From Lemma \ref{le:2.9}, we have $r(F)+2\geq r(H)$, so we have $r(F)=r(H)-2$. Similar to Case 2, we can get that  $C_{p}^{\varphi}$ cannot be of Type E.

Combining with above, we can get that $C_{p}^{\varphi}$ can not be of Type C, D or E, i.e., $p$ is even.

\textbf{Claim 2}. $p\equiv 0(\rm{mod}~4)$.

Suppose to the contrary that $p\equiv 2(\rm{mod}~4)$. Denote by $C_{p}-v=P_{p-1}$, by Lemmas \ref{le:2.7} and \ref{le:2.10}, we have $r(P_{p-1})=r(C_{p})-2$, by Lemma \ref{le:2.8}, we have
$r(G)=r(G-v)+2$.  On the other hand, by (ii) of Lemma \ref{le:4.1}, we have $r(G)=r(G-v)$, a contradiction.

\textbf{Claim 3}. $C_{p}^{\varphi}$ is of Type B.

Suppose to the contrary that $C_{p}^{\varphi}$ is of Type A, by Lemma \ref{le:4.3}, then we have
 $p-2+r(H^{\varphi})=r(\Phi)=r(G)+2\theta(G)$.

 Denote by $C_{p}-v=P_{p-1}$, by Claim 2 and Lemmas \ref{le:2.7}, \ref{le:2.10}, we have $r(P_{p-1})=r(C_{p})$. By Lemma \ref{le:2.8}, we have
$r(G)=r(P_{p-1})+r(H)=p-2+r(H)$.
So, we have $p-2+r(H^{\varphi})=r(\Phi)=r(G)+2\theta(G)=p-2+r(H)+2\theta(H)+2$, i.e., $$r(H^{\varphi})=r(H)+2\theta(H)+2,$$ which contradicts  to the upper bound of $r(H^{\varphi})$ of Theorem \ref{th:3.2}.

\textbf{Claim 4.} Both $H^{\varphi}$ and $F^{\varphi}$ are upper-optimal.

Combining with (ii) of Lemma \ref{le:4.1} and (a) of this theorem, we have
$r(G)+2\theta(G)=r(\Phi)=r(\Phi-v)+2=p+r(F^{\varphi})$, and $r(G)+2\theta(G)=r(G-v)+2\theta(G)=p+r(F)+2\theta(F)$.
So, $r(F^{\varphi})=r(F)+2\theta(F)$, that is, $F^{\varphi}$ is upper-optimal.

Let $u$ be a vertex of $C_{p}$ adjacent to $v$. From (ii) of Lemma \ref{le:4.1}, we have $r(G)=r(G-u)=p-2+r(H)$, $r(G)+2\theta(G)=r(\Phi)=r(\Phi-u)+2=p+r(H^{\varphi})$, that is, $p-2+r(H)+2\theta(H)+2=r(G)+2\theta(G)=p+r(H^{\varphi})$. Then, $r(H^{\varphi})=r(H)+2\theta(H)$, i.e., $H^{\varphi}$ is also upper-optimal.

Thus, (b) of this theorem is obtained.

\textbf{Claim 5.} $r(\Phi)=p+r(F^{\varphi})$, $r(F^{\varphi})=r(H^{\varphi})$, $r(G)=p-2+r(H)$ and $r(F)=r(H)$.

From Claim 4, we can see that $r(\Phi)=p+r(F^{\varphi})=p+r(H^{\varphi})$, that is $r(F^{\varphi})=r(H^{\varphi})$. Note that $r(G)=r(G-u)=r(G-v)$, we have $p-2+r(F)=p-2+r(H)$, i.e., $r(F)=r(H)$. Thus, (c) of this theorem is obtained.

This completes the proof. \quad $\square$


\noindent\begin{theorem}\label{th:4.6}
Let $\Phi=(G, \varphi)$ be a $\mathbb{T}$-gain graph of order $n$. If $\Phi$ is lower-optimal, then
 \begin{enumerate}[(a)]
   \item  Cycles (if any) of $\Phi$ are pairwise vertex-disjoint.
   \item Each cycle (if any) $C_{p}^{\varphi}$ of $\Phi$ is of Type A with order $p\equiv 2(\emph{mod}~4)$.
  \item $r(G)=r(T_{G})+\sum_{O\in\mathcal{C}(G)}(|V(O)|)$ and $r(T_{G})=r(\Gamma_{G})$.
\end{enumerate}
\end{theorem}
\noindent\textbf{Proof.}
If $G$ is an acyclic graph, then the theorem is naturally obtained. Suppose $G$ has cycles, let $v$ be a vertex of some cycle. By (iii) of Lemma \ref{le:4.1}, we know that $v$ lies on only one cycle of $G$, i.e., cycles (if any) of $\Phi$ are pairwise vertex-disjoint.

We now proceed by induction on the order $n$ to prove the left assertions.

If $n=1$, then all left assertions hold naturally. Suppose the left assertions all hold for any lower-optimal $\mathbb{T}$-gain graph order at most $n-1$, and $\Phi$ is a lower-optimal $\mathbb{T}$-gain graph of order $n\geq2$.

\textbf{Case 1.} If $T_{G}$ is an empty graph, i.e., $G$ consists of disjoint cycles and some isolated vertices, then the left assertions follow from the fact that: $\Phi$ is lower-optimal if and only if each component of $\Phi$ is lower-optimal, and a single $\mathbb{T}$-gain cycle $C_{p}^{\varphi}$ is lower-optimal if and only if $C_{p}^{\varphi}$ is of Type A with $p\equiv 2(\textrm{mod}~4)$, as desired.

\textbf{Case 2.} If $T_{G}$ has at least one edge, then $T_{G}$ has at least one pendant vertex $u$. If $u\in U$, then $u$ is also a pendant vertex of $G$, if $u=t_{O}\in C_{G}$, then $G$ has a pendant cycle. In the following, we will consider the following subcases.

\textbf{Subcase 2.1.} $G$ has a pendant vertex $u$.

Let $v$ be the vertex of $G$ adjacent to $u$, $H^{\varphi}=G^{\varphi}-u-v$. By Theorem \ref{th:4.2}, we know that $v$ cannot lie on any cycle of $G$, and $H^{\varphi}$ is lower-optimal. From the induction hypothesis to $H^{\varphi}$, we know that
\begin{enumerate}[(1)]
  \item  Each cycle $C_{p}^{\varphi}$ of $H^{\varphi}$ is of Type A with order $p\equiv 2(\rm{mod}~4)$.
  \item $r(H)=r(T_{H})+\sum_{O\in\mathcal{C}(H)}(|V(O)|)$ and $r(T_{H})=r(\Gamma_{H})$.
\end{enumerate}

Since $H^{\varphi}=G^{\varphi}-u-v$, we can see that all cycles of $G$ belong to $H$, then $\sum_{O\in\mathcal{C}(H)}(|V(O)|)=\sum_{O\in\mathcal{C}(G)}(|V(O)|)$. By (1) of Subcase 2.1, we get each cycle $C_{p}^{\varphi}$ of $\Phi$ is of Type A with order  $p\equiv 2(\rm{mod}~4)$. Noting that $u$ is also a pendant vertex of $T_{G}$ (resp., $\Gamma_{G}$) adjacent to $v$ and $T_{H}=T_{G}-u-v$ (resp., $\Gamma_{H}=\Gamma_{G}-u-v$), then combining with (2) of Subcase 2.1 and Lemma \ref{le:2.7}, we have
$$r(G)=r(H)+2=r(T_{H})+\sum_{O\in\mathcal{C}(H)}(|V(O)|)+2=r(T_{G})+\sum_{O\in\mathcal{C}(G)}(|V(O)|),$$
and
$$r(T_{G})=r(T_{H})+2=r(\Gamma_{H})+2=r(\Gamma_{G}).$$

\textbf{Subcase 2.2.} $G$ has a pendant cycle $C_{p}$.

Let $v$ be the unique vertex of $C_{p}$ of degree 3, $F^{\varphi}=\Phi-C_{p}^{\varphi}$ and $H^{\varphi}=F^{\varphi}+v$. By (b) of  Theorem \ref{th:4.4}, we know that both $H^{\varphi}$ and $F^{\varphi}$ are lower-optimal. The induction hypothesis to $H^{\varphi}$ implies that
\begin{enumerate}[(i)]
  \item  Each cycle $C_{p}^{\varphi}$ of $H^{\varphi}$ is of Type A with order $p\equiv 2(\rm{mod}~4)$.
  \item $r(H)=r(T_{H})+\sum_{O\in\mathcal{C}(H)}(|V(O)|)$.
\end{enumerate}

Combining with (a) of Theorem \ref{th:4.4}, assertion (i) of Subcase 2.2 and the fact that $\mathcal{C}(G)=\mathcal{C}(H)\cup{\{C_{p}\}}$ imply that each cycle of $\Phi$ is of Type A with order  $p\equiv 2(\rm{mod}~4)$. Applying (c) of Theorem \ref{th:4.4} and  (ii) of Subcase 2.2, we have
\begin{equation}\label{13}
 r(G)=p+r(H)=p+r(T_{H})+\sum_{O\in\mathcal{C}(H)}(|V(O)|).
\end{equation}
Since $T_{H}$ is isomorphic to $T_{G}$, i.e., $r(T_{H})=r(T_{G})$ and $p+\sum_{O\in\mathcal{C}(H)}(|V(O))|=\sum_{O\in\mathcal{C}(G)}(|V(O)|)$, we have
\begin{equation}\label{14}
 r(G)=r(T_{G})+\sum_{O\in\mathcal{C}(G)}(|V(O)|),
\end{equation}
which proves the first assertion of (c) of this theorem.

Noting that $\mathcal{C}(G)=\mathcal{C}(F)\cup{\{C_{p}\}}$, then from (c) of Theorem \ref{th:4.4} and Equation (\ref{14}), we have
\begin{equation}\label{15}
 r(T_{G})=r(G)-\sum_{O\in\mathcal{C}(G)}(|V(O)|)=p+r(F)-\sum_{O\in\mathcal{C}(G)}(|V(O)|)=r(F)-\sum_{O\in\mathcal{C}(F)}(|V(O)|).
\end{equation}
Since $F^{\varphi}$ is also lower-optimal, the first assertion of (c) of this theorem applying to $F$ implies that
\begin{equation}\label{16}
 r(T_{F})=r(F)-\sum_{O\in\mathcal{C}(F)}(|V(O)|).
 \end{equation}
Combining with Equations (\ref{15}) and (\ref{16}), we get that
\begin{equation}\label{17}
 r(T_{G})=r(T_{F}).
\end{equation}
The induction hypothesis to $F^{\varphi}$ implies that
\begin{equation}\label{18}
 r(T_{F})=r(\Gamma_{F}).
\end{equation}
Since $\Gamma_{G}=\Gamma_{F}$, i.e., $r(\Gamma_{G})=r(\Gamma_{F})$, combining with equations (\ref{17}) and (\ref{18}), we have $r(T_{G})=r(\Gamma_{G})$.

This completes the proof. \quad $\square$


\noindent\begin{theorem}\label{th:4.7}
Let $\Phi=(G, \varphi)$ be a $\mathbb{T}$-gain graph of order $n$. If $\Phi$ is upper-optimal, then
 \begin{enumerate}[(i)]
  \item  Cycles (if any) of $\Phi$ are pairwise vertex-disjoint.
  \item Each cycle $C_{p}^{\varphi}$ of $\Phi$ is of Type B with order $p\equiv 0(\emph{mod}~4)$.
  \item $r(G)=r(T_{G})+\sum_{O\in\mathcal{C}(G)}(|V(O)|-2)$ and $r(T_{G})=r(\Gamma_{G})$.
\end{enumerate}
\end{theorem}
\noindent\textbf{Proof.}
If $G$ is an acyclic graph, then the theorem is naturally obtained. Suppose that $G$ has cycles, let $v$ be a vertex of some cycle. By (iii) of Lemma \ref{le:4.1}, we know that $v$ lies on only one cycle of $G$, so (i) of this theorem follows.

We now proceed by induction on the order $n$ to prove the left assertions.

If $n=1$, then (ii) and (iii) of this theorem hold naturally. Suppose the left assertions (ii) and (iii) all hold for any upper-optimal $\mathbb{T}$-gain graph order smaller than $n$, and $\Phi$ is a upper-optimal $\mathbb{T}$-gain graph of order $n\geq2$.

\textbf{Case 1.} If $T_{G}$ is an empty graph, i.e., $G$ consists of disjoint cycles and some isolated vertices, then the left assertions follow from the fact that: $\Phi$ is upper-optimal if and only if each component of $\Phi$ is upper-optimal, and a single $\mathbb{T}$-gain cycle $C_{p}^{\varphi}$ is upper-optimal if and only if $C_{p}^{\varphi}$ is of Type B with $p\equiv 0(\textrm{mod}~4)$.

\textbf{Case 2.} If $T_{G}$ has at least one edge, then $T_{G}$ has at least one pendant vertex $u$. If $u\in U$, then $u$ is also a pendant vertex of $G$, if $u=t_{O}\in C_{G}$, then $G$ has a pendant cycle. In the following, we will consider the following subcases.

\textbf{Subcase 2.1.} $G$ has a pendant vertex $u$.

Let $v$ be the vertex of $G$ adjacent to $u$ in $G$, $H^{\varphi}=G^{\varphi}-u-v$. By Theorem \ref{th:4.2}, we know that $v$ cannot lie on any cycle of $G$, and $H^{\varphi}$ is upper-optimal. From the induction hypothesis to $H^{\varphi}$, we know that
\begin{enumerate}[(1)]
  \item  Each cycle $C_{p}^{\varphi}$ of $H^{\varphi}$ is of Type B with order $p\equiv 0(\rm{mod}~4)$.
  \item $r(H)=r(T_{H})+\sum_{O\in\mathcal{C}(H)}(|V(O)|-2)$ and $r(T_{H})=r(\Gamma_{H})$.
\end{enumerate}

Since $H^{\varphi}=G^{\varphi}-u-v$, we can see that all cycles of $G$ belong to $H$, then $\sum_{O\in\mathcal{C}(H)}(|V(O)|-2)=\sum_{O\in\mathcal{C}(G)}(|V(O)|-2)$.  By (1) of Subcase 2.1, we obtain that each cycle $C_{p}^{\varphi}$ of $\Phi$ is of Type B with order  $p\equiv 0(\rm{mod}~4)$. Noting that $u$ is also a pendant vertex of $T_{G}$ (resp., $\Gamma_{G}$) adjacent to $v$ and $T_{H}=T_{G}-u-v$ (resp., $\Gamma_{H}=\Gamma_{G}-u-v$), then combining with (2) of Subcase 2.1 and Lemma \ref{le:2.7}, we have
$$r(G)=r(H)+2=r(T_{H})+\sum_{O\in\mathcal{C}(H)}(|V(O)|-2)+2=r(T_{G})+\sum_{O\in\mathcal{C}(G)}(|V(O)|-2),$$
and
$$r(T_{G})=r(T_{H})+2=r(\Gamma_{H})+2=r(\Gamma_{G}).$$

\textbf{Subcase 2.2.} $G$ has a pendant cycle $C_{p}$.

Let $v$ be the unique vertex of $C_{p}$ in $G$ with $d_{G}(v)=3$, $F^{\varphi}=\Phi-C_{p}^{\varphi}$ and $H^{\varphi}=F^{\varphi}+v$. By (b) of Theorem \ref{th:4.5}, we know that both $H^{\varphi}$ and $F^{\varphi}$ are upper-optimal. The induction hypothesis to $H^{\varphi}$ implies that
\begin{enumerate}[(i)]
  \item  Each cycle $C_{p}^{\varphi}$ of $H^{\varphi}$ is of Type B with order $p\equiv 0(\rm{mod}~4)$.
  \item $r(H)=r(T_{H})+\sum_{O\in\mathcal{C}(H)}(|V(O)|-2)$.
\end{enumerate}

Combining with (a) of Theorem \ref{th:4.5}, (i) of Subcase 2.2 and the fact that $\mathcal{C}(G)=\mathcal{C}(H)\cup{\{C_{p}\}}$, we get that each cycle of $\Phi$ is of Type B with order  $p\equiv 0(\rm{mod}~4)$. Applying (c) of Theorem \ref{th:4.5} and  (ii) of Subcase 2.2, we have
\begin{equation}\label{19}
 r(G)=p-2+r(H)=p-2+r(T_{H})+\sum_{O\in\mathcal{C}(H)}(|V(O)|-2).
\end{equation}
Since $T_{H}$ is isomorphic to $T_{G}$, i.e., $r(T_{H})=r(T_{G})$, and $p-2+\sum_{O\in\mathcal{C}(H)}(|V(O)-2)|=\sum_{O\in\mathcal{C}(G)}(|V(O)|-2)$, we have
\begin{equation}\label{20}
 r(G)=r(T_{G})+\sum_{O\in\mathcal{C}(G)}(|V(O)|-2),
\end{equation}
which implies the first assertion of (iii) of this theorem holds.

Noting that $\mathcal{C}(G)=\mathcal{C}(F)\cup{\{C_{p}\}}$, by (c) of Theorem \ref{th:4.5} and Equation (\ref{20}), then we have
\begin{equation}\label{21}
 r(T_{G})=r(G)-\sum_{O\in\mathcal{C}(G)}(|V(O)|-2)=p-2+r(F)-\sum_{O\in\mathcal{C}(G)}(|V(O)|-2)=r(F)-\sum_{O\in\mathcal{C}(F)}(|V(O)|-2).
\end{equation}
Since $F^{\varphi}$ is also upper-optimal, the first assertion of (iii) of this theorem applying to $F^{\varphi}$ implies that
\begin{equation}\label{22}
 r(T_{F})=r(F)-\sum_{O\in\mathcal{C}(F)}(|V(O)|-2).
\end{equation}
Equations (\ref{21}) and (\ref{22}) implies that
\begin{equation}\label{23}
 r(T_{G})=r(T_{F}).
\end{equation}
The induction hypothesis to $F^{\varphi}$ implies that
\begin{equation}\label{24}
 r(T_{F})=r(\Gamma_{F}).
\end{equation}
Since $\Gamma_{G}=\Gamma_{F}$, i.e., $r(\Gamma_{G})=r(\Gamma_{F})$, combining with Equations (\ref{23}) and (\ref{24}), we have $r(T_{G})=r(\Gamma_{G})$.

This completes the proof. \quad $\square$



\noindent\begin{theorem}\label{th:4.8}
Let $\Phi$ be a finite $\mathbb{T}$-gain graph  without loops and multiple arcs of order $n$. Then $\Phi$ is lower-optimal if and only if all the following conditions  hold:
\begin{enumerate}[(1)]
  \item Cycles (if any) of $\Phi$ are pairwise vertex-disjoint;
  \item Each cycle $C_{p}^{\varphi}$ of $\Phi$ is of Type A with order $p\equiv 2(\emph{mod}~4)$;
  \item A series of $\delta$-transformations can switch $G$ to a crucial subgraph $G_{0}$, where $G_{0}$ is the disjoint union of $\theta(G)$ cycles together with some isolated vertices.
\end{enumerate}
\end{theorem}

\textbf{Sufficiency:}
Suppose that we can get a crucial subgraph of $G_{0}$ after $k$ steps of $\delta$-transformations to $G$, where $G_{0}$ is the disjoint union of $\theta(G)$ cycles together with some isolated vertices. By Lemmas \ref{le:2.2} and \ref{le:2.7}, we have
\begin{equation}\label{25}
 r(\Phi)=2k+r(G_{0}^{\varphi}),~r(G)=2k+r(G_{0}).
\end{equation}

Since each cycle $C_{p}^{\varphi}$ of the crucial subgraph $G_{0}$ of $G$ is of Type A  with order $p\equiv 2(\rm{mod}~4)$, by (a) and (b) of Lemma \ref{le:2.1}, Lemmas \ref{le:2.5} and \ref{le:2.10}, we have
\begin{equation}\label{26}
 r(G_{0}^{\varphi})=\sum_{O\in\mathcal{C}(G)}r(O^{\varphi})=\sum_{O\in\mathcal{C}(G)}r(O)-2\theta(G)=r(G_{0})-2\theta(G).
\end{equation}
By Equations (\ref{25}) and (\ref{26}), we have
$$r(\Phi)=2k+r(G_{0}^{\varphi})=2k+r(G_{0})-2\theta(G)=2k+r(G)-2k-2\theta(G)=r(G)-2\theta(G).$$
That is, $\Phi$ is lower-optimal.

\textbf{Necessity:} Let $\Phi$ be a lower-optimal $\mathbb{T}$-gain graph. By (a) and (b) of Theorem \ref{th:4.6}, we can obtain that (1) and (2) of this theorem all hold. Thus $G$ has precisely $\theta(G)$ cycles, and the acyclic graph $T_{G}$ respect to $G$ is well defined. Now, we will proceed by induction on the order $n$ of $\Phi$ to prove the (3) of this theorem.

\textbf{Case 1.} If $n=1$, then the assertion holds naturally.

\textbf{Case 2.} Suppose the assertion holds for all lower-optimal  graphs with order smaller than $n$, and let $\Phi$ be a lower-optimal  graph of order $n$.

\textbf{Subcase 2.1.} If $T_{G}$ has no edges, then $G$ is the disjoint union of $\theta(G)$ cycles along with some isolated vertices, and the assertion holds naturally.

\textbf{Subcase 2.2.} If  $T_{G}$ has at least one edge, by (c) of Theorem \ref{th:4.6}, we have $$r(T_{G})=r(\Gamma_{G})=r(T_{G}-C_{G}).$$ (b) of Lemma \ref{le:2.11} shows that there is a pendant vertex of $T_{G}$ not in $C_{G}$. Thus $G$ has at least one pendant vertex. Let $u$ be a pendant vertex of $G$ adjacent to a vertex $v$ of $G$. By Theorem \ref{th:4.2}, $v$ does not lie on any cycle of $G$ and the induced subgraph $H^{\varphi}=(H, \varphi)=\Phi-u-v$ of $\Phi$ is also lower-optimal, and also has $\theta(G)$ cycles. The induction hypothesis applying to $H^{\varphi}$ implies that a series of $\delta$-transformations can switch $H$ to a crucial subgraph of $G_{0}$ consisting of $\theta(G)$ disjoint union  cycles together with some isolated vertices. Combining with the first step of transformation applying to $G$ and all the other $\delta$-transformations done latter, we can switch $G$ to the crucial subgraph $G_{0}$.

This completes the proof. \quad $\square$


\noindent\begin{theorem}\label{th:4.9}
Let $\Phi=(G, \varphi)$ be a finite $\mathbb{T}$-gain graph without loops and multiple arcs of order $n$. Then $\Phi$ is upper-optimal if and only if the following conditions all hold:
\begin{enumerate}[(i)]
  \item Cycles (if any) of $\Phi$ are pairwise vertex-disjoint;
  \item Each cycle $C_{p}^{\varphi}$ of $\Phi$ is of Type B with order $p\equiv 0(\emph{mod}~4)$;
  \item A series of $\delta$-transformations can switch $G$ to a crucial subgraph $G_{0}$, where $G_{0}$ is the disjoint union of $\theta(G)$ cycles together with some isolated vertices.
\end{enumerate}
\end{theorem}

Since the process of the proof in Theorem \ref{th:4.9} is similar to the process of the proof in Theorem \ref{th:4.8}, so we omit the proof.

As an application of Theorems \ref{th:3.2}, \ref{th:4.8} and \ref{th:4.9}, we can obtain the following theorem.

\noindent\begin{theorem}\label{th:4.10}
Let $\Phi=(G, \varphi)$ be a finite $\mathbb{T}$-gain graph without loops and multiple arcs of order $n$. Then $$1-\theta(G)\leq\frac{r(\Phi)}{r(G)}\leq1+\theta(G).$$

\begin{enumerate}[(a)]
  \item The equality in the upper bound holds if and only if $\Phi$ is
\begin{enumerate}[(1)]
\item acyclic, or
\item $\Phi$ consists of a $C^{\varphi}_{4}$ of Type  B and some isolated vertices.
\end{enumerate}

\item The equality in the lower bound holds if and only if $\Phi$ is acyclic.
\end{enumerate}
\end{theorem}
\noindent\textbf{Proof.}
Note that for a given finite non-empty graph $G$, i.e., $G$ has at least an edge, then $r(G)\geq2$ ($G$ must contain $K_{2}$ as its
induced subgraph and $r(K_{2})=2$).

By Theorem \ref{th:3.2}, we know that $r(G)-2\theta(G)\leq r(\Phi)\leq r(G)+2\theta(G)$, combining with the fact that $r(G)\geq2$, we have
\begin{equation}\label{27}
 1-\theta(G)\leq 1-\frac{2\theta(G)}{r(G)}\leq\frac{r(\Phi)}{r(G)}\leq 1+\frac{2\theta(G)}{r(G)}\leq 1+\theta(G).
\end{equation}

Next, we will prove (a) holds.

\textbf{Sufficiency:}
If $\Phi$ is acyclic, i.e., $\theta(G)=0$.  By Lemma \ref{le:2.6}, then we have $r(\Phi)=r(G)$, that is $\frac{r(\Phi)}{r(G)}=1=1+\theta(G)$.

If $\Phi$  consists of a $C^{\varphi}_{4}$  of Type  B and some isolated vertices. Combining with Lemmas \ref{le:2.5}, \ref{le:2.10} and the fact that $\theta(G)=1$,  then we have $r(\Phi)=4$, $r(G)=2$. Hence, $\frac{r(\Phi)}{r(G)}=2=1+\theta(G)$.

\textbf{Necessity:} Combining with $\frac{r(\Phi)}{r(G)}= 1+\theta(G)$ and  the right hand of Equation (\ref{27}), we have,
$$1+\theta(G)=\frac{r(\Phi)}{r(G)}\leq 1+\frac{2\theta(G)}{r(G)}\leq 1+\theta(G).$$

The equality holds if and only if $\theta(G)=0$ or $\Phi$ is
upper-optimal, $r(G)=2$ and $\theta(G)\geq1$.

If $\theta(G)=0$, then $\Phi$ is acyclic, as desired.

If $\Phi$ is upper-optimal, $r(G)=2$ and $\theta(G)\geq1$, combining with Lemmas \ref{le:2.5}, \ref{le:2.10}  and Theorem \ref{th:4.9}, then $\Phi$ consists of a $C^{\varphi}_{4}$  of Type  B and some isolated vertices, as desired.

Finally, we will prove (b) holds.

\textbf{Sufficiency:}
The sufficiency is obvious when $\Phi$ is acyclic.

\textbf{Necessity:} Combining with $\frac{r(\Phi)}{r(G)}= 1-\theta(G)$ and  the left hand of Equation (\ref{27}), we have,
$$1-\theta(G)\leq 1-\frac{2\theta(G)}{r(G)}\leq\frac{r(\Phi)}{r(G)}=1-\theta(G).$$

The equality holds if and only if $\theta(G)=0$ or $\Phi$ is
lower-optimal, $r(G)=2$ and $\theta(G)\geq1$.

If $\theta(G)=0$, then $\Phi$ is acyclic, as desired.

If $\Phi$ is
lower-optimal, by the fact of Lemma \ref{le:2.10} and Theorem \ref{th:4.8} and  $r(G)=2$, that is impossible.

This completes the proof. \quad $\square$


\end{document}